\documentstyle[12pt]{article}
\begin{document}

\renewcommand{\thesection}{\arabic{section}.}
\renewcommand{\theequation}{\arabic{equation}}
\renewcommand {\c}  {\'{c}}
\newcommand {\cc} {\v{c}}
\newcommand {\s}  {\v{s}}
\newcommand {\CC} {\v{C}}
\newcommand {\C}  {\'{C}}
\newcommand {\Z}  {\v{Z}}

\baselineskip=24pt
\begin{center}
{\bf {ON PARASTATISTICS DEFINED AS TRIPLE OPERATOR ALGEBRAS}}\\
Stjepan Meljanac $^\dagger$, Marijan Milekovi\c $^\ddagger$ and 
 Marko Stoji\c $^\dagger$ 
\end{center} 

\begin{center}
{\it $^\dagger$  Rudjer Bo\s kovi\c \ Institute , Bijeni\cc ka c.54, 10001 Zagreb,
Croatia}\\
E-mail: meljanac@thphys.irb.hr\\
{\it $^\ddagger$ Prirodoslovno-Matemati\cc ki Fakultet,Zavod za teorijsku
fiziku,\\
Bijeni\cc ka c.32, 10000 Zagreb, Croatia}\\
E-mail: marijan@phy.hr
\end{center}

\begin{center}
{\large \bf Abstract}
\end{center}
\vspace{0.5cm}
\baselineskip=24pt
We unify parastatistics, defined as triple operator algebras represented on 
Fock space, in a simple way using the transition number operators. We express them as
 a normal ordered expansion of creation and annihilation operators. We discuss 
 several examples of parastatistics, particularly Okubo's and Palev's parastatistics connected 
to 
 many-body Wigner quantum systems. We relate them to the notion of extended Haldane 
statistics.\\





\newpage
\begin{center}
{\large \bf Abstract}
\end{center}
\vspace{0.5cm}
\baselineskip=24pt
We unify parastatistics, defined as triple operator algebras represented on 
Fock space, in a simple way using the transition number operators. 
We express them as
 a normal ordered expansion of creation and annihilation operators.
  We discuss 
 several examples of parastatistics, particularly Okubo's and Palev's 
 parastatistics connected to 
 many-body Wigner quantum systems. We relate them to the notion of 
 extended Haldane statistics.


\newpage
\section{Introduction}
Recently, a class of parastatistics (generalizing Bose and Fermi statistics)
 has 
been reformulated in terms of Lie supertriple systems $^1$. Particularly,
Green's parastatistics $^2$ as well as new kinds of parastatistics discovered 
by Palev 
$^{3,15}$ are reproduced. However, in this approach the positive definite Fock space 
representations are not treated.\\
On the other hand, a unified view of all operator algebras represented on Fock 
spaces
has been presented $^4$. Furthermore, the permutation invariant statistics are also 
studied in detail $^5$.\\
 Along the lines of Refs.(4,5), in this paper we unify, in a simple way  triple 
 operator algebras of Ref.(1), represented 
on the Fock spaces, as well as Greenberg's infinite quon statistics $^6$ and 
Govorkov's paraquantization 
 $^7$. Particularly, we present and discuss parastatistics which naturally 
 appears 
in many-body Wigner quantum systems $^3$ and its (bosonic and supersymmetric)
 extension$^{15}$. 
It appears that they are a generalization of Klein-Marshalek algebra $^8$, 
extensively used in 
nuclear physics. We discuss them in the framework of the Haldane's definition of 
statistics $^9$. We point out that none of them is an example of the original
 Haldane exclusion statistics, but can be 
related to the so-called extended Haldane statistics $^{10}$. 
 For each of them we find the extended Haldane statistics parameters .\\


\section{Operator algebra, Fock space realization and statistics}

Let us start with any algebra of M pairs of creation and annihilation operators
$a^{\dagger}_{i}$,  $a_i $, $i=1,2,..M $ ($a^{\dagger}_{i}$
 is Hermitian conjugated to $a_i $ ).
The  algebra is defined by a normally ordered expansion $\Gamma$ 
(generally 
no symmetry principle is assumed)

\begin{equation}
a_ia^{\dagger}_{j}=\Gamma_{ij}(a^{\dagger};a),
\end{equation}
\noindent
with the number operators $N_i$, i.e., $[N_{i},a_{j}^{\dagger}]=
a_{i}^{\dagger}\delta_{ij}$,
$[N_{i},a_{j}]=-a_{i}\delta_{ij}$. In this case no peculiar relations of 
the type 
$a_{i}^m = a_{j}^n $, $i \neq j$ can appear. Then, every monomial in
 $\Gamma_{ij}$, 
Eq.(1), is of the type ($\cdots a^{\dagger}_{j}\cdots a_i \cdots $) and all other 
indices appear in pairs ($\cdots a^{\dagger}_{k}\cdots a_k \cdots $). 
The corresponding coefficients 
of expansion can depend on the total number operator $ N=\sum_{i=1}^M N_i $.\\
We assume that there is a unique vacuum $|0>$ and the corresponding Fock space 
representation. The scalar product  is uniquely defined by $<0|0>=1$ , 
the vacuum 
condition $ a_i|0> = 0 $, $i=1,2,..M, $ and Eq.(1). A general N-particle 
state is a 
linear combination of the vectors ($a^{\dagger}_{i_{1}}\cdots a^{\dagger}_{i_{N}}|0>$), 
$i_{1},\cdots i_{N}=1,2,...M.$ We consider Fock spaces with no state vectors 
of negative 
squared norms. Note that we do not specify any relation between the creation 
(or annihilation) 
operators. They appear as a consequnce of the norm zero vectors (null-vectors) 
in 
Fock space.\\
For fixed N mutually different indices $i_{1},\cdots i_{N}$, we define the $(N! 
\times N!)$ 
hermitian matrix of  scalar products between states 
$(a^{\dagger}_{i_{\pi (1)}} \cdots a^{\dagger}_{i_{\pi (N)}}|0>)$ for all 
permutations 
$ \pi \in S_N $. The number of linearly independent states among them is given 
by 
$d_{i_{1},\cdots i_{N}}={\it rank}{\cal A}(i_{1},\cdots i_{N})$. The set of  
$d_{i_{1},\cdots i_{N}}$  for all possible $i_{1},\cdots i_{N}$ and all integers N  
completely characterizes the statistics and the thermodynamic properties of a ${\it free}$ 
system with the corresponding Fock space $^{10}$.\\ 
If the algebra (1) is permutation invariant $^5$, i.e. 
$\langle \pi \mu|  \pi \nu \rangle=  \langle \mu|\nu \rangle  $, for all $\pi, \mu ,\nu \in 
S_{N} $, 
 all expansion terms in $\Gamma_{ij}$  of the form (symbolically)
$$
\Gamma_{ij}:= \sum ( \underbrace{a^{\dagger}\cdots a^{\dagger}}_{l-r})(a^{\dagger}_{j} 
\underbrace {a^{\dagger}\cdots a^{\dagger}}_{r} \underbrace {a \cdots a }_{s}a_i)
(\underbrace {a \cdots a }_{l-s}) 
$$
\noindent
have the same coefficient
 for all $i,j=1,2,...M$. ( One single relation, for example 
$a_1a^{\dagger}_2 = \Gamma_{12}$, determines the whole algebra.)\\
For the permutation invariant algebras there are several important consequences $^5$.

{\it Consequences}

(i) The matrices ${\cal A}(i_{1},\cdots i_{N})$ and their ranks do not depend on concrete 
indices 
$i_{1},\cdots i_{N}$, but only on the multiplicities $\lambda_i $ of appearance of the same 
indices\\ 
$\lambda_1 \geq\lambda_2\geq ... \geq \lambda_M \geq 0, 
\,|\lambda|= \sum_{i=1}^{M}\lambda_i = N $, i.e. on the partition  $\lambda$ of $N$.

(ii) For mutually different indices $i_{1},\cdots i_{N}$ i.e. $\lambda_1 = \lambda_2 = \cdots  
\lambda_N = 1, $ 
the generic matrix 
${\cal A}_{1^N}$ is
\begin{equation}
{\cal A}_{1^N} = \sum_{\pi \in S_N} c(\pi) R(\pi) ,
\end{equation}
\noindent
where R is the right regular representation of the permutation group $S_N$ and $c(\pi)$ are 
(real) coefficients. In other words, any row (column) of the matrix determines the 
whole matrix ${\cal A}_{1^N}$.

(iii) All matrices ${\cal A}_{\lambda}$ can be simply obtained from ${\cal A}_{1^N}$ $^5$.
To check that the Fock space does not contain states of negative norms, it is sufficient 
to show that only generic matrices are non-negative $^{11}$.

(iv) For permutation invariant algebras there exist the transition number operators 
$N_{ij}$ , $i,j=1,2...M$ with the properties

\begin{equation}
[N_{ij},a^{\dagger}_k]=\delta_{ik}\,a^{\dagger}_j,\quad [N_{ij},a_k]= -\delta_{jk}\,a_i, \quad 
N^{\dagger}_{ij} = N_{ji}, \quad N_{ii} \equiv N_i.
\end{equation}
\noindent
$N_{ij}$ can be presented similarly as $\Gamma_{ij}$, i.e. as a normal ordered expansion
\begin{equation}
N_{ij}=a^{\dagger}_j a_i + \alpha \sum_l a^{\dagger}_la^{\dagger}_j a_ia_l 
+ \beta \sum_l ( a^{\dagger}_la^{\dagger}_j a_la_i + a^{\dagger}_ja^{\dagger}_l a_ia_l ) + 
\gamma \sum_l a^{\dagger}_ja^{\dagger}_l a_la_i + \cdots ,
\end{equation}
\noindent
where $\alpha , \beta ,\gamma  $ are constants which do not depend on the indices $i,j$.

In the next section we show that all permutation invariant statistics considered by 
Okubo $^1$, Palev $^{3,15}$, Greenberg $^6$, Govorkov $^7$ and Klein and Marshalek $^8$
 can be simply 
unified in terms of triple-operator algebras 

\begin{equation}
[[a_i,a^{\dagger}_j]_q,a^{\dagger}_k] = x \,\delta_{ij}\, a^{\dagger}_k + 
y \,\delta_{ik} \,a^{\dagger}_j 
+ z \,\delta_{jk} \, a^{\dagger}_i ,
\end{equation}
\noindent
for all $i,j,k = 1,2,\cdots M$. Here, $x,y,z \in {\bf R}$ are constants, [ , ] denotes 
the commutator and $[a,b]_q = a b - q b a$ is the q-deformed commutator.\\

{\it Remarks}

1. Equation (1), together with the vacuum condition $ a_i|0> = 0 $, uniquely determines 
all  matrices ${\cal A}_{1^N}$ and ${\cal A}_{\lambda}$. However,  equation (1) 
does not imply positive definiteness, which has to be  checked separately.

2. All other triple-operator relations follow from  Eq.(5) via hermiticity of 
creation and annihilation operators, a linear combination of Eq.(5) with indices 
interchanged and, finally,as null-states of matrices ${\cal A}_{\lambda}$.

3. The algebra with the well defined number operators $N_i$ imply that $z=0$ in (5). 
If $z \neq 0$, there exist peculiar relations of the type $a_i^2=a_j^2$ for all 
$i,j=1,2,\cdots M$, although $a^{\dagger}_i|0\rangle $ are linearly independent states. 
Such peculiar algebras are consistent if the Fock space does not contain null-states $^{12}$.

4. We point out that the algebra (5) can be simply written as the normal ordered expansion
\begin{equation}
a_ia^{\dagger}_j = (1 + x N) \delta_{ij} +q a^{\dagger}_ja_i + y N_{ij} +z N_{ji},
\end{equation}
\noindent
where $N_{ij}$ are the transition number operator of form (3,4) and N is the total number 
operator.


\section{Examples}

{\it Example 1.} Green's parastatistics $^2$ can be presented in the form of Eq.(6) 
with $x=z=0 , y=\frac{2}{p}, p \in {\bf N}$ and $q=\pm 1$, i.e.

\begin{equation}
a_ia^{\dagger}_j = \delta_{ij} \mp a^{\dagger}_ja_i  \pm \frac{2}{p}\, N_{ij},
\end{equation}
where the upper (lower) sign corresponds to para-Bose (para-Fermi) statistics.\\
The transition number operator $N_{ij}$ is, up to the second order, given by $^4$ 

\begin{equation}
N_{ij}=a^{\dagger}_ja_i +\frac{p^2}{4(p-1)} \sum_{l}[Y_{jl}]^{\dagger}[Y_{il}] + \cdots ,
\end{equation}
where $Y_{il}=a_ia_l - q\, (\frac{2}{p}-1)\, a_la_i $.

{\it Example 2.} Govorkov's new paraquantization $^7$ is given by $x = z = 0, \\
y=\frac{\lambda}{p}, \lambda = \pm 1 , p \in {\bf N}$ and $q=0$ :

\begin{equation}
a_ia^{\dagger}_j = \delta_{ij} - \frac{\lambda}{p}\, N_{ij} ,
\end{equation}
with the transition number operator, up to the second order,

\begin{equation}
N_{ij}=a^{\dagger}_ja_i +\frac{p^2}{p^2 - \lambda ^2} \sum_{l}
[Y_{jl}]^{\dagger}[Y_{il}] + \cdots ,
\end{equation}
where $Y_{il}=a_ia_l + (\frac{\lambda}{p})\,a_la_i$.

{\it Example 3.} Greenberg's infinite quon statistics $^6$ is given by $x = y = z = 0$ 
and $-1 < q  < 1$ :
\begin{equation}
a_ia^{\dagger}_j = \delta_{ij} + q\, a^{\dagger}_ja_i
\end{equation}
with transition number operator, up to the second order,

\begin{equation}
N_{ij}=a^{\dagger}_ja_i + \frac{1}{1-q^2} \sum_l[Y_{jl}]^{\dagger}[Y_{il}] + \cdots ,
\end{equation}
where $Y_{il}=a_ia_l -q\, a_la_i$.
(The closed form for $N_{ij}$ to all orders and for the general parameter 
$q_{ij}$ is presented in Ref.(13).)

{\it Example 4.} Palev's A statistics (Fermi case), which appears naturally in the 
treatment of many-body Wigner quantum systems $^3$, is described by the following algebra 
$(i,j,k = 1,2,\cdots M)$ :
\begin{equation}
\begin{array}{c}
[\{a_i,a^{\dagger}_j\},a^{\dagger}_k] = \delta_{ik}\,a^{\dagger}_j- \delta_{ij}\,a^{\dagger}_k, 
\\[4mm] 
[\{a_i,a^{\dagger}_j\},a_k]= -\delta_{jk}\,a_i + \delta_{ij}\,a_k ,
\end{array}
\end{equation}
$$ 
\{a_i,a_j\} = \{a^{\dagger}_i,a^{\dagger}_j\} = 0 .
$$
Hereafter, $\{ , \}$ denotes the anticommutator.\\
(In the original algebra, the operators depend on two indices, $a_i \mapsto a_{\alpha i}$, 
but the structure of the algebra depends on the single index. One recovers the 
original algebra with $\delta_{\alpha i,\beta j}= \delta_{\alpha \beta} \delta_{ij}$).
The vacuum conditions are $a_i|0\rangle  = 0$, $ a_ia^{\dagger}_j|0\rangle  =
 p\, \delta_{ij}\,|0\rangle$ for $p \in {\bf N}$. Upon the redefinition of the operators 
 $(a_i,a^{\dagger}_i) \mapsto (\sqrt{p}a_i , \sqrt{p} a^{\dagger}_i)$, we write the above 
algebra 
 as normal ordered expansion with $x = -\frac{1}{p}$, $y=\frac{1}{p}$, $z=0$ and 
 $q=-1$ 
\begin{equation}
a_ia^{\dagger}_j = (1 - \frac{N}{p})\delta_{ij} -  a^{\dagger}_ja_i +  \frac{1}{p}\, N_{ij}.
\end{equation}
The action of the annihilation operators $a_i$ on the Fock states is obtained from the above 
relation,
 Eq.(14). For example,
$$
a_ia^{\dagger}_ja^{\dagger}_k |0\rangle = (1 - \frac{1}{p})(\delta_{ij}\, a^{\dagger}_k - 
\delta_{ik}\, a^{\dagger}_ j)|0\rangle .
$$
It follows that
\begin{equation}
\begin{array}{c}
a_i(a^{\dagger}_i)^2 |0\rangle = 0, \qquad \forall i\\[4mm]
a_ia^{\dagger}_ia^{\dagger}_k |0\rangle =- a_ia^{\dagger}_ka^{\dagger}_i |0\rangle = 
(1 - \frac{1}{p})\,a^{\dagger}_k |0\rangle ,\quad i \neq k.
\end{array}
\end{equation}
\noindent
Hence, we obtain $\{ a_i,a_j \} = \{ a^{\dagger}_i,a^{\dagger}_j \} = 0$.\\
Generally,  for mutually different indices $i_{1},\cdots i_{N}$, we find

\begin{equation}
a_{i_1}a^{\dagger}_{i_1}a^{\dagger}_{i_2}\cdots a^{\dagger}_{i_N} |0\rangle =
(1 -  \frac{N-1}{p})\,a^{\dagger}_{i_2}\cdots a^{\dagger}_{i_N} |0\rangle ,
\end{equation}
\noindent
in accordance with Ref.(3). The Fock space does not contain negative 
norm states if $p \in {\bf N}$.
 The above equation (16) implies that the allowed states are only those with 
 $ N \leq p$, and the states with $ N > p$ are null-states.
 
The transition number operartor $N_{ij}$, up to the second order, is:

\begin{equation}
N_{ij}=a^{\dagger}_ja_i +\frac{1}{(p - 1)} \sum_{l}a^{\dagger}_la^{\dagger}_j
a_ia_l + \frac{2}{(p - 1)(p - 2)}\sum_{l_1,l_2}a^{\dagger}_{l_2}a^{\dagger}_{l_1}
a^{\dagger}_ja_ia_{l_1}a_{l_2} + \cdots
\end{equation}
\noindent
and terminates with ${\it p}$ creation and ${\it p}$ annihilation operator terms.
For example, if $p=2$, the terms with  $(p-2)$ appearing in the denominator 
do not appear at all. The $p \rightarrow \infty $ reproduces the Fermi algebra. 
We note that case $p=1$ reproduces  the Klein-Marshalek algebra $^8$, namely only the 
one-particle states are allowed:

\begin{equation}
a_ia^{\dagger}_j=(1 - N)\, \delta_{ij} , \qquad N=\sum_l a^{\dagger}_la_l.
\end{equation}
In this sense, algebra (14) generalizes the Klein-Marshalek algebra.\\
It is interesting that the Fock space generated by the algebra (14) is equivalent to 
the Fock space generated by the algebra (with the same vacuum condition imposed)

\begin{equation}
a_ia^{\dagger}_j=(1 - \frac{N}{p}) \, (\delta_{ij} - a^{\dagger}_ja_i),
\end{equation}
with the same $N_{ij}$ and $N$ as given by Eq.(17).\\
Furthermore, there are infinitely many algebras leading to  different generic 
matrices, but with the same statistics. They can be represented by 

\begin{equation}
a_ia^{\dagger}_j=f(N) (\delta_{ij} - a^{\dagger}_ja_i),
\end{equation}
with $f(n) >0$, $n < p$ and $f(p)=0$. The simplest choice is the step function 
$f(N)=\Theta(p-N)$  ( $\Theta (x) =0 , \, x \leq 0 $ and $ \Theta (x) = 1 , \, x>0 $ ).

 We point out that the corresponding statistics is Fermi statistics 
 restricted up to $N\leq p$ N-particle states. Hence, the counting rule is simply 
 $D^F(M,N)=\left( \begin{array}{c}
M\\ N
\end{array} \right) $, $N\leq p $ and  $D^F(M,N)=0$ if $N > p $. 
Recall that Haldane $^9$ introduced the statistics parameter ${\it g}$ through the change of the 
single-particle 
${\it Hilbert }$ space dimension ${\it d_n}$ 
$$
g_{n\rightarrow n + \Delta n} = \frac{d_n -d_{n + \Delta n}}{\Delta n },
$$
where ${\it n}$ is the number of particles and ${\it d_n}$ is the dimension of the one-particle 
Hilbert space obtained by keeping the quantum numbers of ${\it (n-1)}$ particles fixed.
In the similar way we define the ${\it extended }$ statistics parameter through the 
change of the available one-particle ${\it Fock }$-subspace dimension $^{10}$.
Therefore, the above 
statistics is characterized by the Haldane statistical parameter $g=1$

\begin{equation}
g_{n\rightarrow n+k}=\frac{d_n -d_{n +k}}{k}=\frac{(M-n+1)-(M-n-k+1)}{k} =1 ,
\end{equation}
if $n+k\leq p$. If $n+k = p+1$, then $g_{n\rightarrow n+k}=\frac{(M-n+1)}{(p-n+1)}$, 
$n=1,2,\cdots p $ is 
fractional but $g$ is not constant any more. Hence, this is not an example for the 
original Haldane statistics for which the statistics parameter is $g = const.$ Moreover, the 
above statistics is also not the statistics of the Karabali-Nair type $^{14}$, where $a_i^p \neq 
0$, $a_i^{p+1}=0$, 
and for any $N\leq M p$ N-particle state is allowed, since from the Eq.(15) we already have 
$a_i^2=0$ and $N\leq p$.

{\it Example 5.} Palev's A statistics$^{15}$ (Bose case) is the counterpart of the algebra (14), 
namely:

\begin{equation}
\begin{array}{c}
[[a_i,a^{\dagger}_j],a^{\dagger}_k] = -\delta_{ik}\,a^{\dagger}_j 
- \delta_{ij}\,a^{\dagger}_k ,\\[4mm] 
[[a_i,a^{\dagger}_j],a_k]= \delta_{jk}\,a_i + \delta_{ij}\,a_k ,
\end{array}
\end{equation}
$$ 
[a_i,a_j] = [a^{\dagger}_i,a^{\dagger}_j] = 0 , \qquad i,j,k = 1,2,\cdots M.
$$
and the vacuum condition $a_ia^{\dagger}_j|0\rangle  =
 p\, \delta_{ij}\,|0\rangle$ . After the redefinition of the operators 
 $(a_i,a^{\dagger}_i) \mapsto (\sqrt{p}a_i , \sqrt{p} a^{\dagger}_i)$, we write the  
 normal ordered expansion of $a_ia^{\dagger}_j$ as \\( $x = y=-\frac{1}{p}$,  $z=0$ ,
 $q=-1$ )

\begin{equation}
a_ia^{\dagger}_j = (1 - \frac{N}{p})\delta_{ij} + a^{\dagger}_ja_i -  \frac{1}{p}\, N_{ij}.
\end{equation}
The action of the annihilation operators $a_i$ on the Fock states is obtained from 
 Eq.(23). For example,
$$
a_ia^{\dagger}_ja^{\dagger}_k |0\rangle = (1 - \frac{1}{p})(\delta_{ij}\, a^{\dagger}_k +
\delta_{ik}\, a^{\dagger}_ j)|0\rangle .
$$
Hence, we obtain $[a_i,a_j] = [a^{\dagger}_i,a^{\dagger}_j] = 0$.\\
Generally, we find 
\begin{equation}
a_i(a^{\dagger}_1)^{n_1}(a^{\dagger}_2)^{n_2}\cdots (a^{\dagger}_M)^{n_M}|0\rangle=
N_i\, (1-\frac{N-1}{p})\,(a^{\dagger}_1)^{n_1}(a^{\dagger}_2)^{n_2} \cdots 
(a^{\dagger}_i)^{n_i-1}
\cdots a^{\dagger}_M)^{n_M}|0\rangle
\end{equation}
where $N=\sum_{i=1}^M n_i$. The Fock space does not contain negative norm states if $p \in {\bf 
N}$.
The above equation (24) implies that the states with $N\leq p$ are allowed and 
the states with $N>p$ are null-states. The transition number operator $N_{ij}$ has the same 
form, with 
the same coefficients as in the Fermi case, Eq.(17), and terminates with p-annihilation and
p-creation operator terms. The limit $p \rightarrow \infty$ reproduces the Bose algebra. 
We note that if $p=1$ the above algebra (23) reproduces the Klein-Marshalek algebra $^8$. Hence, 
this algebra 
 is the Bose generalization of the Klein-Marshalek algebra.\\
There are again infinitely many algebras leading to  different generic 
matrices but of the same ranks, i.e. statistics. They can be represented by 

\begin{equation}
a_ia^{\dagger}_j=f(N) (\delta_{ij} + a^{\dagger}_ja_i),
\end{equation}
with $f(n) >0$, $n < p$ and $f(p)=0$. The simplest choice is the step function mentioned after 
Eq.(20)
 or
$f(N)= 1 - \frac{N}{p}$.
 The corresponding statistics is Bose statistics 
 restricted  to  N-particle states with $N\leq p$. Hence, the counting rule is simply 
 $D^B(M,N)=\left( \begin{array}{c}
M+N-1\\ N
\end{array} \right)$, $N\leq p $ and  $D^B(M,N)=0$ if $N > p $. Therefore, the above 
statistics is characterized by the Haldane statistics parameter $g=0$

\begin{equation}
g_{n\rightarrow n+k}=\frac{d_n -d_{n +k}}{k}=\frac{M - M}{k} =0 ,
\end{equation}
if $n+k \leq p$. If $n+k = p+1$, then $g_{n\rightarrow n+k}=\frac{M}{(p-n+1)}$, 
$n=1,2,\cdots p$, is 
fractional but  not constant . Hence, this is not an example for the 
original Haldane exclusion statistics for which  $g $ should be constant. The 
above statistics is also not  of the Karabali-Nair type $^{14}$, since $a_i^p \neq 0$, 
 $a_i^{p+1}=0$ but $N\leq p$. This would be equivalent only for the single-mode oscillator, 
$M=1$.

{\it Example 6.} The Bose and Fermi restricted algebra of Refs. (1,3) (the super-triple system) 
can be 
defined as

\begin{equation}
[a_I,a^{\dagger}_J]_{q}= (1 -  \frac{N}{p})\delta_{IJ}  -  \frac{(-)^{\sigma (I)\sigma 
(J)}}{p}\, N_{IJ} ,
\end{equation}
$$
 q= (-)^{\sigma (I) \sigma (J)},
$$
$$
 \sigma (I) = \left\{ \begin{array}{ll}
0 & \mbox {if $I=i $ \, (Bose)}\\
1 & \mbox {if $I= \alpha $ \,(Fermi)}
\end{array}
\right. 
$$
\noindent
where the index $ I \doteq   (i = 1,2,\cdots M_B; \alpha = 1,2,\cdots M_F. )$ denotes  bosonic ( 
fermionic) 
oscillator and $N = N_B +N_F$ is the total number operator. \\
Explicitly,
$$
[a_i,a^{\dagger}_j] = (1 - \frac{N}{p})\delta_{ij} - \frac{1}{p} N_{ij},
$$
$$
\{ a_{\alpha},a^{\dagger}_{\beta} \} = (1 - \frac{N}{p})\delta_{\alpha \beta} + 
\frac{1}{p} N_{\alpha \beta}
$$
$$
[a_i,a^{\dagger}_{\alpha}] = - \frac{1}{p} N_{i\alpha} 
$$
$$
[a_{\alpha},a^{\dagger}_i] = - \frac{1}{p}N_{\alpha i}.
$$
The consistency condition for the algebra (27) reads:

\begin{equation}
N_{IJ}a^{\dagger}_K - (-)^{(\sigma (I) + \sigma (J))\sigma (K)}a^{\dagger}_K N_{IJ}= \delta_{IK} 
a^{\dagger}_J .
\end{equation}
For example,
$$
\{ N_{i\alpha} , a_{\gamma} \} = + \delta_{\alpha \gamma} a_i ,
$$
$$
[N_{i\alpha} , a^{\dagger}_j ] = \delta_{ij} a^{\dagger}_{\alpha}.
$$
Notice that $ ( N_{i\alpha} )^2 =0$. Thus, $N_{i\alpha}$ plays the role of supersymmetric 
charge.
Furthermore, it follows that 
$$
[a_i , a_j ] = \{a_{\alpha} , a_{\beta } \} = [a_i , a_{\alpha}] = 0.
$$
The action of the annihilation operators $ a_i ,a_{\alpha} $ on the Fock states is obtained by 
combining 
Eqs.(27) and 
(28). The N-particle states are allowed only if $N\leq p$, with $p$ being an integer.\\
The transition number operators, up to the second order, are basically similar to (17) and read

\begin{equation}
N_{IJ}= a^{\dagger}_Ja_I + \frac{1}{(p - 1)} \sum_{L} (-)^{\sigma(L)(\sigma(I) + \sigma(J))}
a^{\dagger}_La^{\dagger}_Ja_Ia_L +
\end{equation}
$$
+ \frac{2}{(p - 1)(p - 2)}\sum_{L_1,L_2} (-)^{(\sigma(L_1)+ \sigma(L_2))(\sigma(I) + \sigma(J))}
a^{\dagger}_{L_2}a^{\dagger}_{L_1}
a^{\dagger}_Ja_Ia_{L_1}a_{L_2} + \cdots ,
$$
where the sum over $L$ runs over bosonic ($i=1,2,\cdots M_B)$ and fermionic \\$(\alpha 
=1,2,\cdots M_F)$ indices.

In the limit $p \rightarrow \infty$, the above algebra reduces to the ordinary Bose and Fermi 
algebra. 
If $p=1$, the above algebra reduces to the Klein - Marshalek algebra with $M_B + M_F$ 
oscillators.

{\it Example 7.} Okubo's triple operator algebra (Example 4. in Ref.(1)) is defined for the 
fermionic operators $a_i$ as
\begin{equation}
[\{a_i,a^{\dagger}_j \},a^{\dagger}_k]= (\frac{2}{p}) (-\delta_{ij}a^{\dagger}_k 
-\delta_{jk}a^{\dagger}_i 
+\delta_{ik}a^{\dagger}_j).
\end{equation}
The normal ordered expansion of $a_ia^{\dagger}_j$ is given by ($x = z = - \frac{2}{p}$, 
$y=\frac{2}{p}$, $q=-1$)
\begin{equation}
a_ia^{\dagger}_j = (1 - \frac{2 N}{p}) \delta_{ij} - a^{\dagger}_j a_i + (\frac{2}{p}) (N_{ij} - 
N_{ji}).
\end{equation}
In the limit $p \rightarrow \infty$, it becomes the Fermi algebra.\\
From (31) it follows that
\begin{equation}
\begin{array}{c}
a_i(a^{\dagger}_j)^2 |0\rangle = -(\frac{2}{p})a^{\dagger}_i|0\rangle  , \qquad \forall 
i,j\\[4mm]
a_ia^{\dagger}_ia^{\dagger}_k |0\rangle =- a_ia^{\dagger}_ka^{\dagger}_i |0\rangle = 
(1 - \frac{2}{p})\,a^{\dagger}_k |0\rangle \quad i \neq k.
\end{array}
\end{equation}
Therefore, 
\begin{equation}
\begin{array}{c}
\{ a_i,a_j \} = \{ a^{\dagger}_i,a^{\dagger}_j \} = 0 ,\qquad i \neq j, \\[4mm]
(a_i)^2 = A , \qquad [a_i , A] = 0 ,\quad  [a_i , A^{\dagger}]=-(\frac{2}{p})a^{\dagger}_i ,
\quad \forall i ,\\[4mm]
(a_i)^p \neq 0 , \qquad (a_i)^{p+1} = 0.
\end{array}
\end{equation}
However, in the Fock space there are negative norm states since $\langle 0 |(a_i)^2 
(a^{\dagger}_i)^2 | 0 \rangle 
= - (\frac{2}{p}) < 0 $. The necessary condition for absence of such states is $z \geq 0$. The 
algebra similar to 
the algebra described by Eqs.(31-33) but with the positive definite Fock representations has 
been called peculiar 
algebra  and was studied in  Ref.(12).\\
Finally, let us mention that all Lie (super) algebras are triple systems ( since 
$[a_i,a^{\dagger}_j]_{\pm} = 
\delta_{ij} (c_i + d_i N_i)$ )
 and for a irreducible representations characterized with highest (lowest)  
weight state $\Lambda$ ("vacuum") one can find the following normal ordered expansion
$$
a_ia^{\dagger}_i = \Gamma _i (a^{\dagger},a;\Lambda), \qquad a_ia^{\dagger}_j = \pm 
a^{\dagger}_ja_i
$$
However, these systems are not permutation invariant in the sense we defined in this paper.

\bigskip

{\bf Acknowledgement}\\
We thank T.Palev for providing us with references (3,15) and interesting remarks.

\newpage
\baselineskip=24pt
{\bf References}
\begin{description}
\item{1.}
S.Okubo,  J.Math.Phys.{\bf 35}, 2785 (1994); "{\it Super-triple systems and applications to 
parastatistics and Yang-Baxter equations}" (hep-th/9306160).
\item{2.} 
H.S.Green,  Phys.Rev.{\bf 90}, 170 (1953); Y.Ohnuki and S.Kamefuchi, {\it Quantum field theory 
and 
parastatistics}, University of Tokio Press, Tokio, Springer, Berlin, 1982.
\item{3.}
T.Palev, Czech. Journ. Phys.  {\bf B29}, 91 (1979); ibid.  {\bf B32}, 680 (1982);
J.Math.Phys.{\bf 23}, 1778 (1982);  
T.Palev and S.Stoilova , J. Phys.A :Math.Gen.{\bf 27}, 977 (1994);
 ibid. {\bf 27},7387 (1994); 
J.Math.Phys.{\bf 38}, 2806 (1997).
\item{4.}
S.Meljanac and M.Milekovi\c ,  Int.J.Mod.Phys.{\bf A11}, 1391 (1996). 
\item{5.}
B.Meli\c $\,$ and S.Meljanac,  Phys.Lett.{\bf A226 }, 22  (1997); S.Meljanac, M.Stoji\c $\,$ and 
D.Svrtan,  
Phys.Lett.{\bf A224 }, 319  (1997).
\item{6.}
O.W.Greenberg, Phys. Rev.{\bf D43}, 4111 (1991); Phys.Rev.Lett.{\bf 64}, 705 (1990).
\item{7.} 
A.B.Govorkov,  Theor.Math.Phys.{\bf 98}, 107 (1994); Nucl.Phys.{\bf B 365}, 381 (1991).
\item{8.}
A.Klein and E.R.Marshalek, Rev.Mod.Phys.{\bf 63}, 375 (1991); Z.Phys.{\bf A329},
441 (1988). 
\item{9.}
F.D.M.Haldane,  Phys.Rev.Lett.{\bf 67}, 937 (1991).
\item{10.}
S.Meljanac and M.Milekovi\c ,  Mod.Phys.Lett.{\bf A11}, 3081 (1996).
\item{11.}
D.Zagier, Comm.Math.Phys.{\bf 147}, 199 (1992); S.Meljanac and D.Svrtan,  Comm.Math. {\bf 1}, 1 
(1996) 
and preprint IRB-TH-5/95. 
\item{12.}
S.Meljanac, M.Milekovi\c $\,$ and A.Perica,  Europhys.Lett.{\bf 33}, 175 (1996);
S.Meljanac, M.Milekovi\c $\,$ and A.Perica, Int.J.Theor.Phys. {\bf 36}, 11 (1997).
\item{13.}
S.Meljanac and A.Perica, J.Phys.A :Math.Gen.{\bf 27}, 4737 (1994);  Mod.Phys.Lett.{\bf A9}, 3293 
(1994).
\item{14.}
D.Karabali and V.P.Nair, Nucl.Phys.{\bf B 438}, 551 (1995).
\item{15.}
T.Palev,Habilitation Thesis, Sofia 1976; {\it "Lie algebraic aspect of quantum statistics. 
Unitary 
quantization (A-statistics) } ( hep-th/9705032 ).
\end{description}
\end{document}